# A HALL-FUSION BIALGEBRA

Brian J. Day

Macquarie University, NSW 2109, Australia

ABSTRACT. We describe what might be called the "Hall-fusion" bialgebra constructed from a promonoidal double, and mention the corresponding face version for probicategories.

## 1. THE VERTEX VERSION

This was derived from [4] and uses the notion of a promonoidal category [1].

Let $\mathcal{A}$ be a finite (skeletal) $\text{Vect}_k$-category with hom $\mathcal{A}(a, b) = 0$ if $a \neq b$ in $\text{obj}\mathcal{A}$.

Let
$$p : \mathcal{A}^{\text{op}} \otimes \mathcal{A}^{\text{op}} \otimes \mathcal{A} \longrightarrow \text{Vect}_{\text{fd}} \, , \quad I \in \mathcal{A} \, ,$$

$$q : \mathcal{A} \otimes \mathcal{A} \otimes \mathcal{A}^{\text{op}} \longrightarrow \text{Vect}_{\text{fd}} \, , \quad J \in \mathcal{A} \, ,$$

be promonoidal structures on $\mathcal{A}$ and $\mathcal{A}^{\text{op}}$ respectively, and let $(p \otimes q, (I, J))$ be the tensor product promonoidal structure on $\mathcal{A} \otimes \mathcal{A}^{\text{op}}$ [2].

Define the associative algebra $B(p, q)$ to be the $k$-linear span of the set
$$\left\{ e\binom{a}{b} : (a, b) \in \mathcal{A}^{\text{op}} \otimes \mathcal{A} \right\}$$

with product $\mu$ defined by the bilinear extension of
$$e\binom{a}{c} \cdot e\binom{b}{d} = \sum_{u,v} \frac{\dim p(a, b, u) \dim q(c, d, v)}{\dim \mathcal{A}(u, u) \dim \mathcal{A}(v, v)} \, e\binom{u}{v}$$

and unit
$$1 = e\binom{I}{J},$$

and define the coproduct by the linear extension of

(1.1) $$\Delta\left(e\binom{a}{b}\right) = \sum_u e\binom{a}{u} \otimes e\binom{u}{b}$$

with counit

(1.2)
$$\varepsilon\left(e\binom{a}{b}\right) = \delta_{a,b} \ .$$

**Proposition 1.1.** *The above structure defines a bialgebra $B(p, q)$ if*

$$\int^{a,b} p(a, b, u) \otimes q(a, b, v) \cong \mathcal{A}^{\mathrm{op}}(u, v) \ . \ \square$$

Also, given an "antipode" functor $S : \mathcal{A}^{\mathrm{op}} \longrightarrow \mathcal{A}$ with $S^2 = 1$ on $\mathcal{A}$, we obtain an "antipode" $S$ on $B(p, q)$; namely,

$$S\left(e\binom{a}{b}\right) = e\binom{Sb}{Sa} \ ,$$

with $S(x.y) = S y S x$ and $S(1) = 1$ if we define

$$p(a, b, u) = q(S b, S a, S u) \quad \text{and} \quad I = S J \ .$$

The "von Neumann" axiom

$$\mu_3 (1 \otimes S \otimes 1) \Delta_3 = 1 : B(p, q) \longrightarrow B(p, q)$$

is satisfied if also

$$\int^a q_3(S a, a, b; c) \cong \int^a q_3(a, S a, b; c) \cong \mathcal{A}^{\mathrm{op}}(b, c) \ ,$$

where $q_3$ is defined by

$$q_3(a, b, c; d) = \int^u q(a, b, u) \otimes q(u, c, d) \cong \int^u q(b, c, u) \otimes q(a, u, d)$$

and the isomorphism is by associativity of $q$ (see [1]).

This (almost) suffices to make $B(p, q)$ a Hopf algebra. However, the most obvious examples are based on $\mathcal{A}$ being the $k$-linearisation of a finite groupoid, which may here yield nothing more than some familiar ("double"-type) construction on the groupoid algebra over $k$.

## 2. THE FACE VERSION

The face version is derived directly from [4] using (finite) "probicategories" [3] in place of the promonoidal categories of Section 1, the face idempotents being defined by the 1-cell identities (i.e., the 0-cells) of the particular two probicategories under consideration, in much the same manner as that of the "face model" construction in [4] Section 3 p. 233 for a (finite) directed graph.

Explicitly, we consider two (finite) probicategory structures of the form

$$(\mathcal{A}_{ij}, p_{ijk}) \quad \text{and} \quad \left(\mathcal{A}_{ij}^{\mathrm{op}}, q_{ijk}\right)$$

respectively (where the $\mathcal{A}_{ij}$ are suitable $k$-linear categories), both indexed by the same finite set $N = \{i, j, k, \ldots\}$ of 0-cells. We define the face idempotents

$$\mathring{e}_i = \sum_j e\binom{i}{j} \quad \text{and} \quad e_j = \sum_i e\binom{i}{j}$$

for 0-cells $i, j \in N$, the product

$$e\binom{a}{c} \cdot e\binom{b}{d} = \sum_{u,v} \frac{\dim p_{ijk}(a,b,u) \dim q_{ijk}(c,d,v)}{\dim \mathcal{A}_{ik}(u,u) \dim \mathcal{A}_{ik}(v,v)} e\binom{u}{v}$$

for 1-cells $a, c \in \mathcal{A}_{ij}$, $b, d \in \mathcal{A}_{jk}$ (and $u, v \in \mathcal{A}_{ik}$), and also define the coproduct and counit as in (1.1) and (1.2). Then proceed as in Section 1.

**Remark** Clearly, there is no need to have only a finite number of objects in the categories $\mathcal{A}$ (and $\mathcal{A}_{ij}$) introduced above, provided the corresponding promultiplication $p$ (and $p_{ijk}$) etc. have appropriate finite support.